\documentclass[12pt]{amsart} 
\author{Constantin N. Beli}
\title[Reciprocity laws for the Legendre symbols $\leg{a+b\k
m}p$]{Reciprocity laws for Legendre symbols of the type $\leg{a+b\k
m}p$ - long version${}^1$\footnote{${}^1$This paper contains some
rather lenghty discutions involving $2$-adic Hilbert symbols. Readers
who want to skip this technical part may consider the short version.}}  
\usepackage{amsfonts,amsmath,amssymb}
\def\a{\alpha} \def\b{\beta} \def\c{\gamma} 
\def\D{\Delta} \def\d{\delta} \def\e{\varepsilon} 
  \def\h{\frac} 
\def\j{\infty} \def\k{\sqrt}  
\def\m{\lim}    
\def\p{\partial}   
   
 \def\te{\theta}  
  
\def\({\overline} \def\){\underline}
\def\<{\cdot} \def\go{\mathfrak}
\def\>{~~~~~~~} \def\#{{\bf
Definition}} \def\*{\section} \def\be{\begin{equation}}
\def\ee{\end{equation}}

\def\sb{\subset}  \def\sbq{\subseteq}  
\def\ti{\times}   
   
 \def\ff{\dot{F}}  
 \def\mo{{\rm mod}~}  
  \def\fs{\ff^2}  
\def\p{\go p}    
\def\*{\sharp}  \def\0{} 
 \def\1{^{-1}}  
 \def\[{\prec} \def\]{\succ} 
\def\bmat{\left(\begin{array}} \def\emat{\end{array}\right)} \def\ev{\equiv}

 \def\m2{~(\mo 2)} \def\no{\noindent}
 \def\btm{\begin{thm}}
\def\etm{\end{tm}}
 \def\blem{\begin{lem}}
\def\elem{\end{lem}}

\newtheorem{theorem}{Theorem}[section]
\newtheorem{proposition}[theorem]{Proposition}
\newtheorem{lemma}[theorem]{Lemma}
\newtheorem{definition}{Definition}
\newtheorem{corollary}[theorem]{Corollary}

\newtheorem{bof}[theorem]{}
\newtheorem{teorema}{Theorem}

\def\qed{\mbox{$\Box$}\vspace{\baselineskip}}
\def\pf{$Proof.$} 
\def\bco{\begin{corollary}} \def\eco{\end{corollary}} 
\def\bdf{\begin{definition}} \def\edf{\end{definition}} 
\def\btm{\begin{theorem}} \def\etm{\end{theorem}} 
\def\blm{\begin{lemma}} \def\elm{\end{lemma}} 
\def\bff{\begin{bof}\rm} \def\eff{\end{bof}}
\def\btr{\begin{teorema}} \def\etr{\end{teorema}}
\def\bpr{\begin{proposition}} \def\epr{\end{proposition}}

\def\de{\newcommand} \de\tm[1]{{\no\bf Theorem~#1}}
\def\mb{\mathbb} 
\def\RR{{\mb R}}\def\QQ{{\Bbb Q}}\def\ZZ{{\mb
Z}} 
\def\la{\langle} \def\ra{\rangle}

\de\lm[1]{{\no\bf Lemma~#1}}
\de\df[1]{{\no\bf Definition~#1}} \de\co[1]{{\no\bf Corollary~#1}}
\de\tp[1]{\te (#1 )} \de\ts[1]{\te (O^-(#1 ))} \de\ty[1]{\te
(O(#1 ))} \de\tx[1]{\te (#1 )} \de\up[1]{(1+\p^{#1} )\fs}
\de\lr[1]{\longrightarrow^{\!\!\!\!\!\!\!\! #1}}
\de\lf[1]{\longleftarrow^{\!\!\!\!\!\!\!\! #1}}
\de\si[1]{\sim^{\!\!\!\!\! #1}} \de\apr[1]{\approx^{\!\!\!\!\! #1}}
\de\leg[2]{\left(\frac {#1}{#2}\right)} 
 \def\ff{\dot F} \def\fs{\ff^2} \def\ee{\dot E}
 \def\bu{\bullet} \def\ot{\otimes}
\def\xz{\xrightarrow}
\DeclareMathOperator\im{Im} 
 \DeclareMathOperator\ord{ord}

\begin{document}
\maketitle

{\bf Note} This paper is an announcement. We do not prove the main
result, Theorem 2.5 (as well as Theorem 2.10). We only show how
Theorem 2.5 can be used to produce reciprocity laws for Legendre
symbols of the type $\leg{a+b\k m}p$. Note that in the applications we
won't use Theorem 2.10, which is stronger than Theorem 2.5. Therefore
readers who are not interested in details may skip Theorem 2.10 and
everything related to it. (Lemma 1.1 and the part of \S2 following
Theorem 2.5.) 

\section{Introduction and notations}

Given $p>2$ a prime and $a,b,m\in\ZZ/p\ZZ$ such that $\leg
mp=\leg{a^2-mb^2}p=1$, the Legendre's symbol $\leg{a+b\k m}p$ is
defined as $\leg{a+b\a}p$, where $\a\in\ZZ/p\ZZ$ satisfies
$\a^2=m$. The definition is independent of the choice of the square
root $\a$ of $m$ because $(a+b\a )(a-b\a )=a^2-b^2m$ so
$\leg{a+b\a}p\leg{a-b\a}p=\leg{a^2-mb^2}p=1$ so
$\leg{a+b\a}p=\leg{a-b\a}p$. One may aslo define $\leg{a+b\k m}p$ if
$a^2-mb^2=0$ and $a\neq 0$. In this case the two radicals of $m$ are
$\a =\pm\h ab$. If we take $\a =-\h ab$ then $\leg{a+b\a}p=\leg 0p$,
which is not convenient. Hence we will take $\a =\h ab$ and we get
$\leg{a+b\k m}p=\leg{a+b\a}p=\leg{2a}p$. 


In the particular case when $a=0$, $b=1$ we get $\leg{\k m}p$, defined
when $\leg mp=\leg{0^2-m\cdot 1^2}p=1$, i.e. when $p\ev 1\pmod 4$ and
$\leg mp=1$. This coincides with the rational quartic residue symbol
$\leg mp_4$. 

The Legendre symbol $\leg{a+b\k m}p$ can also be defined when $a,b,m$
are $p$-adic integers satisfying $\leg mp=\leg{a^2-mb^2}p=1$ or $p\mid
a^2-mb^2$, $p\nmid a$. 

For any prime $p$ of $\QQ$, including the archimedian prime $p=\j$, we
denote by $\leg{\cdot,\cdot}p:\QQ_p^\ti/(\QQ_p^\ti
)^2\ti\QQ_p^\ti/(\QQ_p^\ti )^2\to\{\pm 1\}$ the Hilbert symbol. (At
$p=\j$ we have $\QQ_\j =\RR$.)

Similarly as for $\leg{a+b\k m}p$ we can introduce the Hilbert symbol
$\leg{a+b\k m,c}p$. It is defined for $a,b,c,m\in\QQ_p$ satisfying
$m\in (\QQ_p^\ti )^2$ and $\leg{a^2-mb^2,c}p=1$. If
$\pm\a\in\QQ_p^\ti$ are the two square roots of $m$ then $\leg{(a+b\a
)(a-b\a ),c}p=\leg{a^2-mb^2,c}p=1$ so $\leg{a+b\a,c}p=\leg{a-b\a ,c}p$
so there is no ambiguity in this definition. If $a^2-mb^2=0$ and
$a\neq 0$, same as for $\leg{a+b\k m}p$, we define $\leg{a+b\k
m,c}p=\leg{2a,c}p$. Throughout this paper we will consider
$\leg{a+b\k m,c}p$ only for $p=2$. 

Various mathematicians have obtained results involving Legendre
symbols of the type $\leg {a+b\k m}p$. Most of this results involve
the biquadratic symbols $\leg mp_4$ and they are called biquadratic or
quartic reciprocity laws. In this paper we will give some very general
reciprocity laws, which generalize all the existing results. In \S3 we
show how many of these results can be obtained as a consequence of
Theorem 2.5, which we only state in this paper.

\bdf For any $A\in\QQ^\ti/(\QQ^\ti )^2$ we denote by $\( A$ the only
squarefree integer such that $A=\( A(\QQ^\ti )^2$,
\edf

\subsection*{The algebras $T(V)$, $S(V)$ and $S'(V)$}${}$

For convenience we denote by $V:=\QQ^\ti/(\QQ^\ti )^2$ and for any
prime $p$ we denote $V_p:=\QQ_p^\ti/(\QQ_p^\ti )^2$. In particular,
$V_\j =\RR^\ti/(\RR^\ti )^2$. 

Now $V$ and $V_p$ are $\ZZ/2\ZZ$-vector spaces. Therefore we may
consider the corresponding tensor algebras $T(V)$ and $T(V_p)$ and the
symmetric algebras $S(V)$ and $S(V_p)$. Note that, while $V$ has a
multiplicative notation, $T(V)$ and $S(V)$ have an additive one. In
order to prevent confusion for any $A_1,\ldots,A_n\in V$ we will
denote by $A_1\bu A_2\bu\cdots\bu A_n$ their product in $S^n(V)$
(instead of simply $A_1\cdots A_n$). For example, the distributivity
law on $S^2(V)$ will be written as $AB\bu C=A\bu C+B\bu C$.

We also define the algebra 
$$S'(V)=T(V)/\la A_1\ot\cdots\ot A_n-A_{\pi (1)}\ot\cdots\ot A_{\pi
(n)}\mid\pi\in A_n\ra.$$
(Same definition as for $S(V)$ but this time the products are
invariant only under even permutations of the factors.) 

We denote by $A_1\odot\cdots\odot A_n$ the image in $S'(V)$ of
$A_1\ot\cdots\ot A_n\in T(V)$. We will only be concerned with the
homogenous component of degree $3$, $S'^3(V)=T^3(V)/\la A\ot B\ot
C-B\ot C\ot A\ra$. We have $A\odot B\odot C=B\odot C\odot A=C\odot
A\odot B$ and $B\odot A\odot C=A\odot C\odot B=C\odot B\odot
A$. However $A\odot B\odot C\ne B\odot A\odot C$ (unless $A,B,C$ are
linearly dependent). Similarly we define $S'(V_p)$. 

For any $\xi\in V$, $T(V)$, $S(V)$ or $S'(V)$ and any prime $p$ we
denote by $\xi_p$ the image of $\xi$ in $V_p$, $T(V_p)$, $S(V_p)$ or
$S'(V_p)$. When there is no danger of confusion we simply write $\xi$
instead of $\xi_p$. 

\blm We have an exact sequence
$$S^3(V)\xz\tau S'^3(V)\xz\rho S^3(V)\to 0,$$
where $\tau$ is given by $A\bu B\bu C\mapsto A\odot B\odot C+B\odot
A\odot C$ and $\rho$ by $A\odot B\odot C\mapsto A\bu B\bu C$.
\elm
\pf The mapping $\rho$ is well defined and surjective because
$S^3(V)=T^3(V)/I$ and $S'^3(V)=T^3(V)/I'$ where $I,I'$ are bubgroups
of $T^3(V)$ with $I'\sb I$. For $\tau$ we note that the mapping
$:V^3\to S'^3(V)$ given by $(A,B,C)\mapsto A\odot B\odot C+B\odot
A\odot C$ is trilinear and also symmetric. (Recall that $A\odot B\odot
C=B\odot C\odot A=C\odot A\odot B$ and $B\odot A\odot C=A\odot C\odot
B=C\odot B\odot A$.) Hence $\tau$ is well defined. 

We have $\rho (\tau (A\bu B\bu C)))=A\bu B\bu C+B\bu A\bu C=2(A\bu
B\bu C)=0$ so $\rho\circ\tau=0$ so $\im\tau\sbq\ker\rho$. For the
reverse inclusion take $\xi =\sum_iA_i\odot B_i\odot
C_i\in\ker\rho$. Then $0=\rho (\xi )=\sum_iA_i\bu B_i\bu C_i$ so there
are $A_{j,1},A_{j,2},A_{j,3}\in V$ and $\pi_j\in S_3$ such that
$\sum_iA_i\ot B_i\ot C_i=\sum_j(A_{j,1}\ot A_{j,2}\ot
A_{j,3}-A_{j,\pi_j(1)}\ot A_{j,\pi_j(2)}\ot A_{j,\pi_j(3)})$. This
implies that $\xi =\sum_iA_i\odot B_i\odot C_i=\sum_j(A_{j,1}\odot
A_{j,2}\odot A_{j,3}-A_{j,\pi_j(1)}\odot A_{j,\pi_j(2)}\odot
A_{j,\pi_j(3)})$. But $A_{j,\pi_j(1)}\odot A_{j,\pi_j(2)}\odot
A_{j,\pi_j(3)}$ equals $A_{j,1}\odot A_{j,2}\odot A_{j,3}$ or
$A_{j,2}\odot A_{j,1}\odot A_{j,3}$ if $\pi_j$ is even or odd,
respectively. So $A_{j,1}\odot A_{j,2}\odot
A_{j,3}-A_{j,\pi_j(1)}\odot A_{j,\pi_j(2)}\odot A_{j,\pi_j(3)}$ is $0$
if $\pi_j\in A_3$ and it is $A_{j,1}\odot A_{j,2}\odot
A_{j,3}-A_{j,2}\odot A_{j,1}\odot A_{j,3}=\tau (A_{j,1}\bu A_{j,2}\bu
A_{j,3})$ otherwise. It follows that $\xi =\tau (\eta )$, where $\eta
=\sum_{j,\pi_j\notin A_3}A_{j,1}\bu A_{j,2}\bu A_{j,3}$. \qed

\section{Main result}

\bdf We denote by $\mathcal D$ the set of all $(B,A,C)\in V^3$ such that

1) $\leg{A,B}p=\leg{A,C}p=\leg{B,C}p=1$ for every prime $p$, including
the archimedian prime $p=\j$.  

2) $(\( A,\( B,\( C)=1$.

3) At least one of $\( A,\( B,\( C$ is $\ev 1\pmod 4$. 
\edf

\bdf We define $f_1:{\mathcal D}\to\{\pm 1\}$ as follows. If
$(B,A,C)\in\mathcal D$ we take $x,y,z\in\ZZ$ with $(x,y,z)=1$ such
that $x^2-\( Ay^2=\( Bz^2$. (The existence of such $x,y,z$ is ensured
by Minkovski's theorem since $\leg{A,B}p=1$ for all $p$.) Then define
$f_1(B,A,C)=\a_\j\prod_{p\mid 2\( C}\a_p$, where
$$\a_\j =\begin{cases}1&\text{if }C>0\\
sgn(x)&\text{if }C<0\end{cases},$$

$$\a_p=\begin{cases}\leg{x+y\k{\( A}}p&\text{if }p\nmid\( A\\
\leg{2(x+z\k{\( B})}p&\text{if }p\nmid\( B\end{cases}$$
if $p\mid\( C$, $p>2$ and

$$\a_2=\begin{cases}1&\text{if }\( C\ev 1\pmod 8\\
\leg{x+y\k{\( A},C}2&\text{if }\( A\ev 1\pmod 8\\
\leg{2(x+z\k{\( B}),C}2&\text{if }\( B\ev 1\pmod 8\\
\leg{z,5}2&\text{if }\( A\ev\( C\ev 5\pmod 8\\
-\leg{y,5}2&\text{if }\( B\ev\( C\ev 5\pmod 8\\
-\leg{3x+z\k{5\( B},-1}2&\text{if }\( A\ev\( B\ev 5\pmod 8,~\( C\ev
-1\pmod 8\\
-\leg{5x+z\k{5\( B},3}2&\text{if }\( A\ev\( B\ev 5\pmod 8,~\( C\ev
3\pmod 8\
\end{cases}.$$
\edf

\bdf We define $f_2:{\mathcal D}\to\{\pm 1\}$ as follows. If
$(B,A,C)\in\mathcal D$ we take $x,y,z\in\ZZ$ with $(x,y,z)=1$ such
that $x^2-\( By^2=-\({ABC}z^2$. (The existence of such $x,y,z$ is
ensured by Minkovski's theorem since $\leg{A,B}p=\leg{C,B}p=1$ so
$\leg{-ABC,B}p=1$ for all $p$.) Then define
$f_2(B,A,C)=\b_\j\prod_{p\mid 2\( C}\b_p$, where
$$\b_\j =\begin{cases}1&\text{if }C>0\\
sgn(x)&\text{if }C<0\end{cases},$$

$$\b_p=\begin{cases}\leg xp&\text{if }p\nmid\( A\\
\leg{2(x+y\k{\( B})}p&\text{if }p\nmid\( B
\end{cases}$$
if $p\mid\( C$, $p>2$ and
$$\b_2=\begin{cases}1&\text{if }\( C\ev 1\pmod 8\\
\leg{x,C}2&\text{if }\( A\ev 1\pmod 8\\
\leg{2(x+y\k{\( B}),C}2&\text{if }\( B\ev 1\pmod 8\\
\leg{y+z\k{\h{\({ABC}}{\( B}},5}2&\text{if }\( A\ev\( C\ev 5\pmod 8\\
-\leg{z,5}2&\text{if }\( B\ev\( C\ev 5\pmod 8\\
\leg{-2(5x+y\k{5\( B}),C}2&\text{if }\( A\ev\( B\ev 5\pmod 8
\end{cases}.$$
\edf
\bff{\bf Note} If $x=0$ and $\( A\ev 1\pmod 8$ then we make the
convention that in the definition of $\b_2$ we don't use the option
$\b_2=\leg{x,C}2=\leg{0,C}2$, which is not defined. More precisely, if
$x=0$ then $0^2-\( By^2=-\({ABC}z^2$, which implies that $B=ABC$ so
$A=C$. Hence if $\( A\ev 1\pmod 8$ then also $\( C\ev 1\pmod 8$ so we
will choose the option $\b_2=1$. (See the special case 2.)
\eff

\bff{\bf Remark} In the definition of $f_1$ and $f_2$ we may replace
the condition that $x,y,z$ are relatively prime integers by the
condition that $x,y,z$ are relatively prime elements of $\ZZ [\h
12]$. Indeed, if we replace $x,y,z$ by $2^kx,2^ky,2^kz$ for some
$k\in\ZZ$ then $\a_\j$ is not changed, $\a_2$ is replaced by
$\leg{2^k,C}2\a_2$ and for any $p>2$, $p\mid\( C$ $\a_p$ is replced by
$\leg{2^k}p\a_p$. Since $\leg{2^k,C}2\prod_{p\mid\(
C,p>2}\leg{2^k}p=\prod_p\leg{2^k,C}p=1$ the product $\a_\j\prod_{p\mid
2\( C}\a_p$ is not changed. Similarly for the product
$\b_\j\prod_{p\mid 2\( C}\b_p$, which defines $f_2(B,A,C)$. 
\eff

\bff{\bf Approximations of 2-adic square roots} In the definition of
$\a_2,\b_2$ for $\( C\not\ev 1\pmod 8$ we have formulas of the type
$\leg{a+b\k m,C}2$, where $a,b\in\QQ_2$ and $m\in 1+8\ZZ_2$. Here we
have the liberty of choosing either of the two quadratic roots of
$m$. Let $m=1+8x$. We make the convention that $\k m$ is the quare
root of $m$ that is $\ev 1\pmod 4$. We have $\k m+(1-4x)\ev 2\pmod 4$
and $(\k m)^2-(1-4x)^2=-16x(x-1)\vdots 32$ so $\k
m-(1-4x)\vdots 16$ so $\k m\ev 1-4x=\h{3-m}2\pmod{16}$. 

Consider the $\pm$ sign such that 
$$\ord_2(a\pm b\cdot\h{3-m}2)=
\min\{\ord_2(a+b\cdot\h{3-m}2),\ord_2(a-b\cdot\h{3-m}2)\}.$$
It follows
that $\ord_2(a\pm b\cdot\h{3-m}2)\leq\ord_2((a+b\cdot\h{3-m}2)+
(a-b\cdot\h{3-m}2))=\ord_22b\cdot\h{3-m}2=\ord_22b$. (Recall that
$\h{3-m}2=1-4x$ is an odd 2-adic integer.) Hence $a\pm
b\cdot\h{3-m}2\mid 2b$. Since also $b(\k m-\h{3-m}2)\vdots 16b$ we get
$$\h{a\pm b\k m}{a\pm b\cdot\h{3-m}2}-1=\h{\pm b(\k
m-\h{3-m}2)}{a\pm b\cdot\h{3-m}2}\vdots \h{16b}{2b}=8.$$
It follows that $\h{a\pm b\k m}{a\pm b\cdot\h{3-m}2}$ is a square so
$\leg{a\pm b\k m,C}2=\leg{a\pm b\cdot\h{3-m}2,C}2$. 

If $\( C\ev 5\pmod 8$ then our formula becomes $\leg{a\pm b\k
m,5}2$. By the same reasoning as above, since $\k m\ev 1\pmod 4$, if
we take the $\pm$ sign such that $\ord_2(a\pm
b)=\min\{\ord_2(a+b),\ord_2(a-b)\}$ then $\h {a\pm b\k m}{a\pm b}\ev
1\pmod 2$ so $\leg{\h {a\pm b\k m}{a\pm b},5}2=1$, which implies that
$\leg{a\pm b\k m,5}2=\leg{a\pm b,5}2$.
\eff

Although calculating the Hasse symbols from the definition of
$\a_2,\b_2$ can be done very easily in numerical cases by using
2.3, proving general results may be quite labourios. Readers who want
to skip this part may check the short version of this paper.

{\bf Special cases}

{\bf 1.} $B=-A$ and we want to calculate $f_1(-A,A,C)$. Condition 1)
from the definition of $\mathcal D$ means that $\leg{A,C}p=\leg{-1,C}p=1$
$\forall p$ and condition 2) means that $(\( A,\( C)=1$. We have $\(
B=-\( A$ and we can take $x=0$, $y=z=1$. Since $\leg{-1,C}\j =1$ we
have $C>0$ so $\a_\j =1$. If $p\mid\( C$, $p>2$ then $\a_p=\leg{0+\k{\(
A}}p=\leg{\( A}p_4$. Also 
$$\a_2=\begin{cases}1&\text{if }\( C\ev 1\pmod 8\\
\leg{\k{\( A},C}2&\text{if }\( A\ev 1\pmod 8\\
\leg{2\k{-\( A},C}2&\text{if }\( A\ev -1\pmod 8\\
1&\text{if }\( A\ev\( C\ev 5\pmod 8\\
-1&\text{if }\( A\ev 3\pmod 8,~\( C\ev 5\pmod 8\\
\end{cases}.$$
(The proof follows straight-forward from the definition of
$\a_2$. Since $\( B=-\( A$ the cases with $\( A\ev\( B\ev 5\pmod 8$ do
not occur. The condition $\( B\ev 1\pmod 8$ means $\( A\ev -1\pmod 8$
and the condition $\( B\ev\( C\ev 5\pmod 8$ means $\( A\ev 3\pmod 8$,
$\( C\ev 5\pmod 8$.)

{\bf 2.} $A=C$ and we want to calculate $f_2(B,C,C)$. Condition 1) in
the definition of $\mathcal D$ means $\leg{B,C}p=\leg{-1,C}p=1$ $\forall
p$ and condition 2) means that $(\( B,\( C)=1$. We have
$\({ABC}=\({BC^2}=\( B$ so we are looking for $x,y,z$ relatively
prime, such that $x^2-\( By^2=-\(B z^2$. One obvious choice is
$(x,y,z)=(0,1,1)$ but it is more convenable to choose $(x,y,z)=(0,\h
12,\h 12)$. (See Remark 2.2.) Since $\leg{-1,C}\j =1$ we have $C>0$ so
$\b_\j=1$. For $p\mid\( C$, $p>2$ we have $p\nmid\( B$ so we take
$\b_p=\leg{2(x+y\k{\( B})}p=\leg{2(0+\h 12\k{\( B})}p=\leg{\k{\(
B}}p=\leg{\( B}p_4$. If $\( A=\( C$ is odd then $\leg{C,-1}2=1$
implies $\( C\ev 1\pmod 4$. If $\( C\ev 1\pmod 8$ then $\b_2=1$. If
$\( A=\( C\ev 5\pmod 8$ then $\b_2=\leg{y+z\k{\h{\({ABC}}{\(
B}},5}2=\leg{\h 12+\h 12\k{\h{\( B}{\( B}},5}2=\leg{1,5}2=1$. If $\(
A=\( C$ is even then by the definition of $\mathcal D$ we have $\( B\ev
1\pmod 4$, which, together with $\leg{B,C}2=1$, implies $\( B\ev
1\pmod 8$. (We have $2\|\( C$.) It follows that $\b_2=\leg{2(x+y\k{\(
B}),C}2=\leg{2(0+\h 12\k{\( B}),C}2=\leg{\k{\( B},C}2$. Note that
$\k{\( B}$ is an odd $2$-adic integer so by multiplying with $\pm 1$
we may assume that $\k{\( B}\ev 1\pmod 4$. Hence $\b_2=\leg{\k{\(
B},C}2$ is $1$ or $-1$ according as $\k{\( B}\ev 1$ or $5\pmod
8$. But by 2.3 $\k{\( B}\ev\h{3-\( B}2\pmod 8$ so $\b_2=\leg{\k{\(
B},C}2$ is $1$ if $\( B\ev 1\pmod{16}$ and it is $-1$ if $\( B\ev
9\pmod{16}$. For short, if $\( C$ is even then $\b_2=\leg{\( 
B}2_4$, where $\leg\cdot 2_4:1+8\ZZ_2\to\{\pm 1\}$ is given by $\leg
a2_4=1$ if $a\ev 1\pmod{16}$ and $\leg a2_4=-1$ if $a\ev
9\pmod{16}$. If $\( C$ is odd then $\b_2=1$. In conclusion:
$$f_2(B,C,C)=\prod_{p\mid\( C}\leg{\( B}p_4=:\leg{\( B}{\( C}_4.$$

{\bf 3.} $ABC=-1$, i.e. $A=-BC$, and we want to calculate
$f_2(B,-BC,C)$. In this case condition 1) in the definition of $\mathcal
D$ is equivalent to $\leg{B,C}p=1$ $\forall p$, while condition 2) is
vacuous. We have $\({ABC}=1$ so we need $x,y,z$ such that $x^2-\(
By^2=z^2$. An obvious choice is $x=z=1$, $y=0$. Since
$x>0$ we have $\b_\j =1$. If $p\mid\( C$, $p>2$
then $p\nmid\( A=-\({BC}$ is equivalent to $p\mid\( B$. In this case
$\b_p=\leg xp=1$. If $p\nmid\( B$ then $\b_2=\leg{2(x+y\k{\( B})}p=\leg
2p$. Also
$$\b_2=\begin{cases}\leg{2,C}2&\text{if }\( B\ev 1\pmod 8\\
\leg{-2,C}2&\text{if }\( B\ev 5\pmod 8\\
1&\text{otherwise }\end{cases}.$$
(Since $x=z=1$, $y=0$ and $\({ABC}=-1$ we get $\b_2=1$ if $\( C\ev
1\pmod 8$ or $\( A\ev 1\pmod 8$; $\b_2=\leg{2,C}2$ if $\( B\ev 1\pmod
8$; $\b_2=\leg{\k{\h{-1}{\( B}},5}2=1$ if $\( A\ev\( C\ev 5\pmod 8$
(in this case $\( B\ev -1\pmod 8$ so $\h{-1}{\( B}$ is a unit in $\ZZ_2$
and so is $\k{\h{-1}{\( B}}$); $\b_2=-1$ if $\( B\ev\( C\ev 5\pmod 8$;
$\b_2=\leg{-10,-1}2=-1$ if $\( A\ev\( B\ev 5\pmod 8$ (in this case
$\( C\ev -1\pmod 8$). Note that if $\( B\ev 5\pmod 8$ from
$\leg{B,C}2=1$ we get that $\( C$ is odd. If $\( B\ev 5,\( C\ev 1\pmod
8$ then $\b_2=1$; if $\( B\ev 5,\( C\ev 5\pmod 8$ then $\b_2=-1$; if
$\( B\ev 5,\( C\ev -1\pmod 8$, which is equivalent to $\( B\ev\( A\ev
5\pmod 8$, then $\b_2=-1$;  if  $\( B\ev
5,\( C\ev 3\pmod 8$, which is equivalent to $\( B\ev 5\pmod 8$, $\(
A\ev 1\pmod 8$, then $\b_2=1$. Hence if $\( B\ev 5\pmod 8$ then
$\b_2=\leg{-2,C}2$. If $\( B\ev 1\pmod 8$ then $\b_2=\leg{2,C}2$. In
all the other cases $\b_2=1$.)

Note that $\prod_{p\mid\( C,p>2}\b_p=\prod_{p\mid\( C,p\nmid 2\(B}\leg
2p=(-1)^{\h{n^2-1}8}$, where $n=\h{\( C}{(\( C,2\( B)}$. Since also
$\b_\j =1$ we get $f_2(B,-BC,C)=(-1)^{\h{n^2-1}8}\b_2$. 

\bff{\bf Remark} If $p\mid\( C$, $p>2$ and $p\nmid\( A,\( B$ then in
the definition of $\a_p$ and $\b_p$ we can choose the formula from
both the $p\nmid\( A$ case and the $p\nmid\( B$ case. However the
outcome is the same. 

Indeed, for $\a_p$ we have $x^2-\( Ay^2=\( Bz^2$ so $(x+y\k{\(
A})\cdot 2(x+z\k{\( B})=(x+y\k{\( A}+z\k{\ B})^2$. It follows that
$\leg{(x+y\k{\( A})\cdot 2(x+z\k{\( B})}p=1$ so $\leg{x+y\k{\(
A}}p=\leg{2(x+z\k{\( B})}p$. 

For $\b_p$ we have $x^2-\(
By^2=\({ABC}z^2$ so $x\cdot 2(x+y\k{\( B})=(x+y\k{\(
B})^2-\({ABC}z^2\ev (x+y\k{\( B})^2\pmod p$. It follows that
$\leg{x\cdot 2(x+y\k{\( B})}p=1$ so $\leg xp=\leg{2(x+y\k{\( B})}p$. 

We have similar redundancies in the definition of $\a_2,\b_2$ but
again all the cases of the definition which apply produce the same
outcome. 
\eff

We now state our main result. 

\btm (i) The functions $f_1$ and $f_2$ are well defined, i.e. they are
independent of the choice of $x,y,z$, and they are equal. We denote
$f=f_1=f_2$. 

(ii) $f$ is symmetric in all three variables. 

(iii) If $(B_i,A_i,C_i)\in\mathcal D$ for $1\leq i\leq n$ and
$\sum_iB_i\odot A_i\odot C_i=0$ then $\prod_if(B_i,A_i,C_i)=1$. 
\etm

Theorem 2.5 justifies the following definition.

\bdf We denote by $W$ the subspace of $S'^3(V)$ generated by $B\odot
A\odot C$ with $(B,A,C)\in\mathcal D$. We define the group morphism $\c
:W\to\{\pm 1\}$ by $B\odot A\odot C\mapsto f(B,A,C)$. 
\edf

\bff{\bf Remark} Part of Theorem 2.5(ii) follows from Theorem
2.5(iii). Namely, since $B\odot A\odot C=A\odot C\odot B=C\odot B\odot
A$, we have $f(B,A,C)=f(A,C,B)=f(C,B,A)$, i.e. $f$ has a circular
symmetry. However from Theorem 2.5(ii) we know that $f$ is symmetric,
not merely circular symmetric. Hence one may assume that we have a
more precise result, namely $\prod_if(B_i,A_i,C_i)=1$ whenever
$\sum_iB_i\bu A_i\bu C_i=0$ in $S^3(V)$. The reason that this
doesn't happen is the following. Assume that $(B_i,A_i,C_i)\in\mathcal D$
and $\sum_iB_i\bu A_i\bu C_i=0$. Then $\xi :=\sum_i B_i\odot A_i\odot
C_i$ can be written as $\xi =\sum_j(B_j'\odot A_j'\odot C_j'+A_j'\odot
B_j'\odot C_j')$ for some $A_j',B_j',C_j'\in V$. We have 
$\prod_if(B_i,A_i,C_i)=\c (\xi )$. If $(B_j',A_j',C_j')\in\mathcal D$ for
all $j$ then $\c (\xi
)=\prod_jf(B_j',A_j',C_j')f(A_j',B_j',C_j')$. But from Theorem 2.5(ii)
we have $f(B_j',A_j',C_j')=f(A_j',B_j',C_j')$. Hence we get $\c (\xi
)=1$, as expected. The reason why this ``proof'' is wrong is that not
always $B_j',A_j',C_j'$ can be chosen such that
$(B_j',A_j',C_j')\in\mathcal D$. 
\eff

\bff{\bf Example} Assume that $p,q,r,s,a,b,c,d$ are odd primes with
$p\ev q\ev r\ev s\ev 1\pmod 4$, $\leg ap=\leg bp=\leg cp=\leg dp=\leg
aq=\leg dq=\leg br=\leg dr=\leg cs=\leg ds=-1$ and $\leg bq=\leg
cq=\leg ar=\leg cr=\leg as=\leg bs=1$. Then
$(ad,pq,pq),(bd,pr,pr),(cd,ps,ps),(a,rs,rs),(b,qs,qs),\\ (c,qr,qr),
(abcd,pqrs,pqrs)\in\mathcal D$. Indeed, in all cases the condition 2)
from the definition of $\mathcal D$ is trivial and the conditon 3)
follows from $p\ev q\ev r\ev s\ev 1\pmod 4$. Condition 1) for the
first triplet means $\leg{ad,pq}t=1$ and $\leg{pq,-1}t=1$ for any
prime $t$. The second condition follows from the fact that $pq$ is a
sum of two squares. The first condition at $t=\j$ follows from $pq>0$
and at $t=2$ from the fact that $pq\ev 1\pmod 4$ and $ad$ is odd. At
$t=p,q,a$ and $d$ it means $\leg{ad}p=\leg{ad}q=\leg{pq}a=\leg{pq}d=1$
and it follows from $\leg ap=\leg pa=\leg dp=\leg pd=\leg aq=\leg
qa=\leg dq=\leg qd=-1$. For $t\neq\j,2,p,q,a,d$ our statement is
trivial. Similarly for the other triplets. 

We also have $ad\bu pq\bu pq+bd\bu pr\bu pr+cd\bu ps\bu ps+a\bu rs\bu
rs+b\bu qs\bu qs+c\bu qr\bu qr+abcd\bu pqrs\bu pqrs=0$. However the
product
$f(ad,pq,pq)f(bd,pr,pr)f(cd,ps,ps)f(a,rs,rs)f(b,qs,qs)f(c,qr,qr)\\
f(abcd,pqrs,pqrs)$ is $-1$. Indeed, by the special case 2 our product
is equal to $\leg{ad}{pq}_4\leg{bd}{pr}_4\leg{cd}{ps}_4\leg
a{rs}_4\leg b{qs}_4\leg c{qr}_4\leg{abcd}{pqrs}_4=\\
\leg{a^2b^2c^2d^4}p_4\leg{a^2b^2c^2d^2}q_4\leg{a^2b^2c^2d^2}r_4
\leg{a^2b^2c^2d^2}s_4=\\ \leg{abc}p\leg{abcd}q\leg{abcd}r\leg{abcd}s
=(-1)\cdot 1\cdot 1\cdot 1=-1$. 
\eff 


\blm If $U=\tau\1 (W)$ then for every $p$ the image $U_p$ of $U\sbq
S^3(V)$ in $S^3(V_p)$ is generated by elements of the form $A\bu A\bu
B$ with $A,B\in V_p$. 
\elm

\bff{\bf Remark} With the exception of the case $p=2$ we have
$\dim_{\ZZ/2\ZZ}V_p\leq 2$ so $U_p=S^3(V_p)$. If $p>2$ then
$V_p=\la\D_p,p\ra$, where $\D_p$ is a nonsquare unit in $\ZZ_p$ so
$S^3(V_p)=\la\D_p\bu\D_p\bu\D_p,\D_p\bu\D_p\bu p,p\bu p\bu\D_p,p\bu
p\bu p\ra =U_p$. If $p=\j$ then $V_\j =\RR^\ti/(\RR^\ti )^2=\la -1\ra$
so $S^3(V_\j )=\la (-1)\bu (-1)\bu (-1)\ra =U_\j$. 

If $p=2$ then $V_2=\la -2,3,6\ra$ so a basis for $S^3(V_2)$ is $\{
A\bu B\bu C\mid A,B,C\in\{ -2,3,6\}\}$. A basis of $U_2$ is made of
all the elements of the basis for $S^3(V_2)$ except $(-2)\bu 3\bu 6$
so $U_2\sb S^2(V_2)$. (The dimensions of $S^3(V_2)$ and $U_2$ are $10$
and $9$, respectively.) 
\eff

\bdf We define $\d :U\to\{\pm 1\}$ by $\d (\xi )=\prod_p\d_p(\xi_p)$,
where the product is taken over all primes, including $p=\j$, and $\d_p
:U_p\to\{\pm 1\}$ is given by $A\bu A\bu B\mapsto\leg{A,B}p$. 
\edf

\btm $\c\circ\tau =\d$.
\etm

Theorem 2.10 provides a generalization of Theorem 2.5(iii). Namely, if
$(B_i,A_i,C_i)\in\mathcal D$ and $\sum_iB_i\bu A_i\bu C_i=0$ then we write
$\xi =\sum_iB_i\odot A_i\odot C_i\in W$ and we have $\rho (\xi
)=0$. It follows that $\xi\in\ker\rho =\im\tau$. (See Lemma 1.1.) Then
$\xi =\tau (\eta )$ for some $\eta\in\tau\1 (W)=U$ and we have
$\prod_if(B_i,A_i,C_i)=\c (\xi )=\c (\tau (\eta ))=\d (\eta )$. Note
that, in principle, calculating $\d (\eta )$ is easier than $\c (\xi
)$ as it only involves usual Legendre symbols. In order to calculate
$\c (\xi )$ one has to compute $f(B_i,A_i,C_i)$, which involves
finding nontrivial zeros for some ternary quadratic form and
calculating Legendre symbols of the type $\leg{a+b\k m}p$. 

Some more explicit formulas for $\d_p$, so for $\d$ are given bellow:

\bff If $p=\j$ then $\d_\j :U_\j\to\{\pm 1\}$ is given by
$$\d_\j (A\bu B\bu C)=\begin{cases}-1&\text{if }A,B,C<0\\
1&\text{otherwise}\end{cases}.$$

If $p>2$ then $\d_p:U_p\to\{\pm 1\}$ is obtained as follows. For any
$A,B,C\in V$ if we write $\( A=p^ra$, $\( B=p^sb$ and $\( C=p^tc$,
where $r,s,t\in\{ 0,1\}$ and $p\nmid abc$ then 
$$\begin{array}{ll}\d_p(A\bu B\bu C)=&\leg{-1}p^{rst}\leg ap^{st}\leg
bp^{rt}\leg cp^{rs}\cdot\\ {}&\cdot\left(\leg bp*\leg cp\right)^r
\left(\leg cp*\leg ap\right)^s\left(\leg ap*\leg
bp\right)^t,\end{array}$$ 
where $*:\{\pm 1\}\ti\{\pm 1\}\to\{\pm 1\}$ is given by
$$\e *\eta =\begin{cases}-1&\text{if }\e =\eta =-1\\
1&\text{otherwise}\end{cases}.$$

If $p=2$ and $\eta\in U$ then $\eta_2=\sum_iA_i\bu B_i\bu C_i$ with
$A_i,B_i,C_i\in V_2$. We write $A_i,B_i,C_i$ in terms of the basis
$-2,3,6$ of $V_2$. We have $A_i=(-2)^{r_{i,1}}3^{r_{i,2}}6^{r_{i,3}}$,
$B_i=(-2)^{s_{i,1}}3^{s_{i,2}}6^{s_{i,3}}$ and
$C_i=(-2)^{t_{i,1}}3^{t_{i,2}}6^{t_{i,3}}$ with
$r_{i,j},s_{i,j},t_{i,j}\in\{ 0,1\}$. Then
$$\d_2 (\eta )=\prod_i\prod_{j=1}^3(-1)^{r_{i,j}s_{i,j}t_{i,j}}.$$

Note that we may extend $\d_2:U_2\to\{\pm 1\}$ to the whole $S^3(V_2)$
by setting arbitrarily $\d_2 ((-2)\bu 3\bu 6)=1$. This way $\d
:U\to\{\pm 1\}$ extends to the whole $S^3(V)$. The formula above for
$\d_2$ will hold for for all $\eta\in S^3(V)$. If
$A=(-2)^{r_1}3^{r_2}6^{r_3}$, $B=(-2)^{s_1}3^{s_2}6^{s_3}$ and
$C=(-2)^{t_1}3^{t_2}6^{t_3}$ then $\d_2 (A\bu B\bu
C)=\prod_{j=1}^3(-1)^{r_js_jt_j}$.  
\eff

In most applications we won't need Theorem 2.10. In fact in all proofs
from the next section we don't even need the full strength of Theorem
2.5(iii). It is enough to use a weaker version of Theorem 2.5(iii)
where the condition $\sum_iB_i\odot A_i\odot C_i=0$ is replaced by
$\sum_iB_i\ot A_i\ot C_i=0$. (Note that $\sum_iB_i\bu A_i\bu
C_i=0\Longrightarrow\sum_iB_i\odot A_i\odot
C_i=0\Longrightarrow\sum_iB_i\ot A_i\ot C_i=0$.) The only symmetry
properties we need are the ones following from Theorem 2.5(ii).

\section{Applications}

We now recover some results regarding quartic reciprocity that can be
found in [L, \S5] or on wikipedia at

http://en.wikipedia.org/wiki/Quartic$\_\,$reciprocity



\subsection*{Formulas for $\leg mp_4$}

By the special case 2 if $p\ev 1\pmod 4$ is a prime and $m=\pm
q_1\cdots q_k$ with $\leg{q_i}p=1$ then $(m,p,p)\in\mathcal D$ and
$$f_2(m,p,p)=\leg mp_4.$$ 

By Theorem 2.5(i) and (ii) we have $\leg
mp_4=f_2(m,p,p)=f_1(p,p,m)$. We have $p=a^2+b^2$ with $2\mid b$. Then
$p^2-b^2p=a^2p$ so in the definition of $f_1(p,p,m)$ we may take
$(x,y,z)=(p,b,a)$. We have $f_1(p,p,m)=\a_\j\prod_{q\mid 2m}\a_q$. Now
$x=p>0$ so $\a_\j =1$ and if $q\mid m$, $q>2$ then $q\nmid p=\( A$ so
$\a_q=\leg{x+y\k{\( A}}q=\leg{p+b\k p}q$.

Assume first that $m=2$. Then $\leg 2p=1$ so $p\ev 1\pmod 8$,
i.e. $4\mid b$. We have $\leg 2p_4=f_1(p,p,2)=\a_\j\a_2=\a_2$. Since
$\( A=p\ev 1\pmod 8$ we have $\a_2=\leg{x+y\k{\( A}}2=\leg{p+b\k
p,2}2$. Now $\k p$ is an odd integer in $\ZZ_2$ and so $4\mid b$ so
$b\k p\ev b\pmod 8$ so $p+b\k p\ev p+b\ev 1+b\pmod 8$. If $b\ev 0\pmod
8$ then $p+b\k p$ is $\ev 1\pmod 8$ so it is a square in $\ZZ_2^\ti$
so $\leg{p+b\k p,2}2=1$ so $\leg 2p_4=1$. If $b\ev 4\pmod 8$ then
$p+b\k p\ev 1+4=5\pmod 8$ so $\leg{p+b\k p,2}2=\leg{5,2}2=-1$ so $\leg
2p_4=-1$. In conclusion $\leg 2p_4=1$ iff $8\mid b$, which is one of
Euler's conjectures, proved by Gauss in 1828.

Take now $m=q^*:=(-1)^{\h{q-1}2}q$, where $q>2$ is prime, $\leg
qp=1$. Then $\leg{q^*}p_4=f_1(p,p,q^*)=\a_\j\a_q\a_2=\a_q\a_2$. We
have $\a_q=\leg{p+a\k p}q$ and since $\( A=p$ and $\( C=q^*$ and $p\ev
q^*\ev 1\pmod 4$ we have three cases for $\a_2$. If $q^*\ev 1\pmod 8$
then $\a_2=1$. If $p\ev 1\pmod 8$ then $\a_2=\leg{x+y\k{\(
A},C}2=\leg{p+b\k p,q^*}2$. But $p$ is odd and $b$ is even so $p+b\k p$
is an odd $2$-adic integer. Since also $q^*\ev 1\pmod 4$ we get again
that $\a_2=1$. If $p\ev q^*\ev 5\pmod 8$ then
$\a_2=\leg{x,5}2=\leg{p,5}2=1$. So $\a_2=1$ and
$\leg{q^*}p=\a_q=\leg{p+b\k p}q$. This result belongs to Lehmer. (See
[L, \S5.4, p. 167].)

\subsection*{Burde}

Assume that $p\ev q\ev 1\pmod 4$ and $\leg pq =1$ and write
$p=a^2+b^2$, $q=c^2+d^2$ with $b$ and $d$ even. Then $pq=e^2+f^2$,
where $e=ac-bd$, $f=ad+bc$ and $e$ is odd and $f$ even. 

Now $p\odot p\odot q+p\odot q\odot q=p\odot pq\odot q$ so by Theorem
2.5(iii)
$$\leg pq_4\leg qp_4=f(p,p,q)f(q,q,p)=f(p,pq,q).$$
Now $f(p,pq,q)=f_2(pq,p,q)$. In the definition of $f_2$ we have $A=p$,
$B=pq$ and $C=q$ so $ABC=1$. Since $e^2-pq=-f^2$ we may take
$(x,y,z)=(e,1,f)$. We have $C=q>0$ so $\b_\j =1$ and the only primes
dividing $2\( C=2q$ are $2,q$ so $f_2(pq,p,q)=\b_q\b_2$. Now $q\nmid
p=\( A$ so $\b_q=\leg xq=\leg eq$. We have $\( A=pq$, $\( C=q$ so $\(
A\ev\( C\ev 1\pmod 4$. We have theree cases. If $\( C\ev 1\pmod 8$
then $\b_2=1$. If $\( C\ev 1\pmod 8$ then $\b_2=\leg{x,C}2=1$. ($x=e$
is odd and $\( C\ev 1\pmod 4$.) If $\( A\ev\( C\ev 5\pmod 8$ then
$\b_2=\leg{y+z\k{\h{\({ABC}}{\( B}},5}2=\leg{y\pm z,5}2$. (See 2.3.)
Here the $\pm$ sign is taken such that $\ord_2 (y\pm z)$ is
minimum. But for both choices of the $\pm$ sign $y\pm z=1\pm f$ is an
odd integer so again $\b_2=1$. Hence $\leg pq_4\leg qp_4=\b_q\b_2=\leg
eq\cdot 1=\leg{ac-bd}q$, which is Burde's law. (See [L, \S5.4,
p. 167].)

We now prove the other formulas for $\leg pq_4\leg qp_4$ from [L,
Theorem 5.7].

We write $\leg pq_4\leg qp_4=f(p,pq,q)=f_2(p,pq,q)$. We have $\(
A=pq$, $\( B=p$, and $\( C=q$ so $\({ABC}=1$. Since $a^2-p=-b^2$ we
can take $(x,y,z)=(a,1,b)$. Again $f_2(p,pq,q)=\b_q\b_2$. Since
$q\nmid\( B=p$ we have $\b_q=\leg{2(x+y\k{\( B})}q=\leg{\h 12(a+\k
p)}q$ and, by the same proof as in the previous case, $\b_2=1$. Hence
$\leg pq_4\leg qp_4=\leg{\h 12(a+\k p)}q$. To obtain $\leg pq_4\leg
qp_4=f_2(p,pq,q)=\leg{b+\k p}q$ we write $b^2-p=-a^2$ and we have
$(x,y,z)=(b,1,a)$. Same as before, $\b_q=\leg{\h 12(b+\k p)}q$, but
this time $\b_2=\leg 2q$, i.e. $\b_2=1$ if $\( C=q\ev 1\pmod 8$ and
$\b_2=-1$ if $\( C\ev 5\pmod 8$. Alternatively we use the fact that in
$\ZZ_q$ $2(a+\k p)(b+\k p)=(a+b+\k p)^2$ so $\leg{\h 12(a+\k
p)}q=\leg{b+\k p}q$. Similarly $\leg pq_4\leg qp_4=\leg{\h 12(c+\k
q)}p=\leg{d+\k p}q$.

In order to obtain $\leg pq_4\leg
qp_4=f(p,pq,q)=\leg{a+b\k{-1}}q=\leg{c+d\k{-1}}p$ we note that $p\odot
pq\odot q=p\odot (-pq)\odot q+p\odot (-1)\odot q$ so
$$\leg pq_4\leg qp_4=f(p,pq,q)=f(p,-pq,q)f(p,-1,q).$$
Now $f(p,-pq,q)=f_2(-pq,p,q)$ and since $A=p$, $B=-pq$ and $C=q$ we
have $ABC=-1$ so we can use the special case 3. We have $n=\h{\(
C}{(\( C,2\( B)}=1$ and, since $\( B=-pq\not\ev 1\pmod 4$,
$\b_2=1$. It follows that $f_2(-pq,p,q)=(-1)^{\h{n^2-1}8}\b_2=1$. So
$$f(p,pq,q)=f(p,-1,q)=f_1(p,-1,q).$$
We have $A=-1$, $B=p$ and $C=q$. Since $a^2+b^2=p$ we may take
$(x,y,z)=(a,b,p)$. We have $C=q>0$ so $\a_\j =1$ so
$f(p,-1,q)=\a_q\a_2$. Since $q\nmid -1=\( A$ we have
$\a_q=\leg{a+b\k{-1}}q$. We have $\( B=p$ an $\( C=q$ ao $\(
B\ev\( C\ev 1\pmod 4$. There are three cases. If $\( C\ev 1\pmod 8$
then $\a_2=1$. If $\( B\ev 1\pmod 8$ then $\a_2=\leg{2(x+z\k{\(
B}),C}2=\leg{2(a+\k p),C}2$. But $\( C\ev 1\pmod 4$ so by 2.3 we have
$\a_2=\leg{2(a\pm 1),C}2$, where the $\pm$ sign is chosen such that
$\ord_2(a\pm 1)$ is minimum. But $a$ is odd so we must have $a\pm 1\ev
2\pmod 4$ (so that $\ord_2(a\pm 1)=1$). Then $\h{a\pm 1}2$ is odd and
since $\( C\ev 1\pmod 4$ we have $\a_2=\leg{2(a\pm 1),C}2=\leg{\h{a\pm
1}2,C}2=1$. If $\( B\ev\( C\ev 5\pmod 8$ then
$\a_2=-\leg{y,5}2=-\leg{b,5}2$. But $a^2+b^2=p=\( C\ev 5\pmod 8$ and
$b$ is even so $b\ev 2\pmod 4$. It follows that $\leg{b,5}2=-1$ and
$\a_2=1$. Since $\a_2=1$ in all cases
$f_1(p,-1,q)=\a_q=\leg{a+b\k{-1}}q$. Note that since $q\nmid p=\( B$
we may also take $\a_q=\leg{2(a+\k p)}q$ and thus we recover the
equality $\leg pq_4\leg qp_4=\leg{\h 12(a+\k p)}q$ from
above. (Alernatively we may use the relation $2(a+b\k{-1})(a+\k
p)=(a+b\k{-1}+\k p)^2$ from Remark 2.4, which implies that $\leg{a+b\k
{-1}}q=\leg{2(a+\k p)}q$ so the two statements are equivalent.) 

To prove [L, Ex. 5.5, p. 176] we note that we also may write
$$\leg pq_4\leg qp_4=f(p,-1,q)=f_1(p,q,-1).$$
We have $A=q$, $B=p$ and $C=-1$ so $f_1(p,q,-1)=\a_\j\a_2$. If
$e^2=pf^2+qg^2$ then $e^2-qg^2=pf^2$ so we may take
$(x,y,z)=(e,g,f)$. We want to prove that
$f_1(p,q,-1)=(-1)^{\h{fg}2}\leg{-1}e$. By permuting, if necessary $p$
and $q$ we may assume that if $p\ev q\ev 5\pmod 8$ then $f$ is even
and $g$ is odd and if $p\not\ev q\pmod 8$ then $p\ev 5,q\ev 1\pmod
8$. We have $\( A=q$, $\( B=p$ so $\( A\ev\( B\ev 1\pmod 4$ and $\(
C=-1$. Since $C<0$ we have $\a_\j =sgn(e)$. There are two cases:
$q\ev 1\pmod 8$ and $q\ev p\ev 5\pmod 8$. In the second case by our
assumption $g$ is even. In the first case $\( A\ev 1\pmod 8$ so
$\a_2=\leg{x+y\k{\( A},C}2=\leg{e+g\k q,-1}2$. If $g$ is even then
$\a_2=\leg{e,-1}2\leg{1+\h ge\k q,-1}2$. Now $\h 1e\k q$ is an odd
$2$-adic integer and $g$ is even. If $2\| g$ then $1+\h ge\k q\ev
3\pmod 4$ so $\leg{1+\h ge\k q,-1}2=-1$. If $4\mid g$ then $1+\h ge\k
q\ev 3\pmod 4$ so $\leg{1+\h ge\k q,-1}2=1$. Hence $\leg{1+\h ge\k
q,-1}2=(-1)^{\h g2}=(-1)^{\h{fg}2}$ and so
$\a_2=(-1)^{\h{fg}2}\leg{e,-1}2$. If $g$ is odd so $f$ is even note
that $pf^2=e^2-qf^2\ev 1-1=0\pmod 8$ so $4\mid f$. Now $e$ and $g\k q$
are odd $2$-adic integers so by replacing, if necessary $\k q$ by $-\k
q$ we have $e\ev g\k q\pmod 4$. It follows that $2\| e+g\k q$. Since
also $(e-g\k q)(e+g\k q)=e^2-qg^2=pf^2\vdots 16$ we get $e-g\k q\vdots
8$ so $e+g\k q\ev 2e\pmod 8$. It follows that $\h{e+g\k q}{2e}\ev
1\pmod 4$ so $\leg{\h{e+g\k q}{2e},-1}2=1$, which implies that
$\a_2=\leg{e+g\k
q,-1}2=\leg{2e,-1}2=\leg{e,-1}2=(-1)^{\h{fg}2}\leg{e,-1}2$. (Recall,
$4\mid f$.) In the second case $\( A\ev\( B\ev 5\pmod 8$ and $\( C=-1$
so $\a_2=-\leg{3x+z\k{5\(
B},-1}2=-\leg{3e+f\k{5p},-1}2=-\leg{3e,-1}2\leg{1+\h f{3e}\k
q,-1}2$. Now $pf^2=e^2-qg^2\ev 1-5\ev 4\pmod 8$ so $2\| f$. Since also
$\h 1{3e}\k{5p}$ is an odd $2$-adic integer we have $1+\h
f{3e}\k{5p}\ev 3\pmod 4$ so $\leg{1+\h f{3e}\k q,-1}2=-1$. Since also
$-\leg{3e,-1}2=\leg{e,-1}2$ we get
$\a_2=-\leg{e,-1}2=(-1)^{\h{fg}2}\leg{e,-1}2$. In conclusion $\leg
pq_4\leg qp_4=\a_\j\a_2=sgn(e)\cdot
(-1)^{\h{fg}2}\leg{e,-1}2=(-1)^{\h{fg}2}\leg{-1}e$. 

Next we prove that if $p=r^2+qs^2$ then $\leg pq_4\leg qp_4=\leg
2q^s$ (see [L, Ex. 5.6, p. 176]). This time we use
$$\leg pq_4\leg qp_4=f(p,-1,q)=f_2(p,q,-1).$$
We have $A=q$, $B=p$ and $C=-1$ so $ABC=-pq$. Since $p^2-pr^2=pqs^2$
we may take $(x,y,z)=(p,r,s)$. We have $f_2(p,q,-1)=\b_\j\b_2$. Since
$x=p>0$ we have $\b_\j =1$ so we have to prove that $\b_2=\leg
2q^s$. Since $C=-1$ there are three cases for $\b_2$. If $q=\( A\ev
1\pmod 8$ then $\b_2=\leg{x,C}2=\leg{p,-1}2=1=\leg 2q^s$. If $q\ev
5\pmod 8$ and $p=\( B\ev 1\pmod 8$ then $\leg 2q^s=(-1)^s$ and
$\b_2=\leg{2(x+y\k{\( B}),C}2=\leg{2(p+r\k p),-1}2=\leg{p+r\k
p,-1}2$. If $s$ is odd and $r$ is even then $r^2=qs^2-p\ev 1-5\ev
4\pmod 8$ so $2\| r$. We have $\b_2=\leg{p,-1}2\leg{1+r\k{\h
1p},-1}2$. But $2\| r$ and $\k{\h 1p}$ is an odd integer so $1+r\k{\h
1p}\ev 3\pmod 4$ so $\leg{1+r\k{\h 1p},-1}2=-1$. It follows that
$\b_2=\leg{p,-1}2\cdot (-1)=-1=\leg 2q^s$. (Recall that $s$ is odd and
$q\ev 5\pmod 8$.) If $s$ is even and $r$ odd then $qs^2=p-r^2\ev
1-1=0\pmod 8$ so $4\mid s$. Since $p$, $r\k p$ are $2$-adic odd
integers, by multiplying $\k p$ with $\pm 1$ we may assume that $p+r\k
p\ev 2\pmod 4$, i.e. $2\| p+r\k p$. Together with
$p^2-pr^2=pqs^2\vdots 16$ we get that $8\mid p-r\k p$ so $p+r\k p\ev
2p\pmod 8$ so $\h 1{2p}(p+r\k p)\ev 1\pmod 4$ so $\leg{\h 1{2p}(p+r\k
p),-1}2=1$, which implies that $\b_2=\leg{2p,-1}2=1=\leg 2s^s$ ($s$ is
even). Now assume that $p\ev q\ev 5\pmod 8$, i.e. $\( A\ev\( B\ev
5\pmod 8$, so $\b_2=\leg{-2(5x+y\k{5\(
B}),C}2=\leg{-2(5p+r\k{5p}),-1}2=-\leg{5p+r\k{5p},-1}2$. If $s$ is 
odd and $r$ even then $r^2=qs^2-p\ev 5-5=0\pmod 8$ so $4\mid r$ so
$5p+r\k{5p}\ev 5p\ev 1\pmod 4$ so $\b_2=-\leg{5p+r\k{5p},-1}2=-1=\leg
2q^s$. (Recall that $s$ is odd and $q\ev 5\pmod 8$.) If
$s$ is even and $r$ is odd then $qs^2=p-r^2\ev 5-1=4\pmod 8$ so $2\|
s$. By multiplying $\k{5p}$ with $\pm 1$ we may assume that $2\|
5p+r\k{5p}$. We also have $25p^2-5pr^2=20p^2+5pqs^2\ev 4+4\ev
8\pmod{16}$. ($p^2\ev 1\pmod 4$ so $20p^2\ev 20\ev 4\pmod{16}$ and
$5pq(\h s2)^2\ev 1\pmod 4$ so $5pqs^2\ev 4\pmod{16}$.) Hence
$5p-r\k{5p}\ev 4\pmod 8$ so $5p+r\k{5p}\ev 10p+4\ev 50+4\ev 6\pmod 8$,
which implies that $\leg{5p+r\k{5p},-1}2=-1$ so $\b_2=1=\leg 2q^s$
($s$ is even).

\subsection*{Scholz}
Assume that $p\ev q\ev 1\pmod 4$ are primes such that $\leg pq=1$ and
assume that $\e$ is a unit of $\QQ (\k p)$ of norm $-1$. Then Scholz's
law states that $\leg\e q=\leg pq_4\leg qp_4$. (See [L, p. 167].) We
could prove Scholz's law directly if $\e\in\ZZ [\leg{1+\k
p}2]\setminus\ZZ [\k p]$. However in order to overcome the $2$-adic
complications we replace $\e$ by $\e^3$ and we may assume that $\e
=t+u\k p\in\ZZ [\k p]$. We have
$$\leg pq_4\leg qp_4=f(p,-1,q)=f_1(-1,p,q)$$
so $A=p$, $B=-1$, $C=q$. Since $t^2-pu^2=-1$ we may take
$(x,y,z)=(t,u,1)$. We have $f_1(-1,p,q)=\a_\j\a_q\a_2$. But $C>0$ so
$\a_\j =1$ and $q\nmid p=\( A$ so $\a_p=\leg{x+y\k{\( A}}q=\leg{t+u\k
p}p=\leg\e p$. So we have to prove that $\a_2=1$. If $\( C=q\ev 1\pmod
8$ then $\a_2=1$ by definition. If $\( C=q\ev 5\pmod 8$ and $\( A
=p\ev 1\pmod 8$ then $\a_2=\leg{x+y\k {\( A},C}2=\leg{t+u\k
p,5}2$. But $t$ and $u$ have opposite parities so $t+u\k p$ is an odd
$2$-adic integer, which implies that $\a_2=1$. If $p\ev q\ev 5\pmod
8$, i.e. $\( A\ev\( C\ev 5\pmod 8$ then
$\a_2=\leg{z,5}2=\leg{1,5}2=1$. 
\vskip 3mm

A similar result holds if $q=2$. Namely, if $p\ev 1\pmod 8$ then
$\leg{1+\k 2}p=\leg 2p_4\leg p2_4$. (See [L, p. 169].) We repeat the
reasoning used to prove that $\leg pq_4\leg qp_4=f(p,-1,q)$. By the
special case 2 we have $f_2(2,p,p)=\leg 2p_2$ and $f_2(p,2,2)=\leg
p2_4$ so $\leg 2p_4\leg p2_4=f(2,2p,p)$. We have
$f(2,2p,p)=f(2,-2p,p)f(2,-1,p)$. But $f(2,-2p,p)=f_2(2,-2p,p)$ can be
calculated using the special case 3 with $B=2$, $C=p$. We have
$n=\leg{\( C}{\( C,2\( B}=p$ and $\( B=2\not\ev 1\pmod 4$ so
$\b_2=1$. Thus $f_2(2,-2p,p)=(-1)^{\h{p^2-1}8}\cdot 1=1$ so $\leg
2p_4\leg p2_4=f(2,2p,p)=f(2,-1,p)=f_1(-1,2,p)$. We have $A=2$, $B=-1$,
$C=p$ so $f_1(-1,2,p)=\a_\j\a_p\a_2$. We have $1^2-2=-1$ so we may
take $(x,y,z)=(1,1,1)$. Since $C=p>0$ we have $\a_\j =1$ and since
$p\nmid 2=\( A$ we may take $\a_p=\leg{x+y\k{\( A}}p=\leg{1
+\k 2}p$ so we have to prove that $\a_2=1$. But this follows from $\(
C=p\ev 1\pmod 8$. 

\section*{References}

[L] Lemmermeyer, Franz, {\it Reciprocity Laws: from Euler to
Eisenstein}, Berlin: Springer-Verlag, 2000

%
%
%
%
%
%
%
%
%
%
%
%
%
%
%
%
%
%
%
%
%
%
%
%
%
%
%

\end{document}